\documentclass{amsart}
\usepackage{amsthm}
\usepackage{amscd}
\usepackage[final]{graphics}
\usepackage{amssymb,amsmath}
\usepackage{latexsym}
\usepackage[english]{babel}\usepackage{tabularx}
\usepackage{xypic}
\usepackage{longtable}

\newtheorem{lem}{Lemma}
    \newtheorem{prop}[lem]{Proposition}
    \newtheorem{thm}[lem]{Theorem}
    \newtheorem{cor}[lem]{Corollary}
 
\theoremstyle{remark}
    \newtheorem{rem}[lem]{Remark}

\DeclareMathOperator{\Sym}{Sym}
\DeclareMathOperator{\im}{Im}

\newcommand{\Sp}{\mathop{\mathrm {Sp}}\nolimits}
\newcommand{\GL}{\mathop{\mathrm {GL}}\nolimits}
\def\s{\mathfrak S}

\newcommand{\op}{\operatorname}
\def\cA{\mathcal A}
\def\cX{\mathcal X}
\def\polybeta{\gamma}

\newcommand{\ab}[1]{\cA_{{#1}}}
\newcommand{\Vor}{{\op{Vor}}}
\newcommand{\Sat}{{\op{Sat}}}
\newcommand{\AVOR}[1]{{\cA_{#1}^{\Vor}}}
\newcommand{\XVOR}[1]{{\cX_{#1}^{\Vor}}}

\newcommand{\ASAT}[1]{{\cA_{#1}^{\Sat}}}
\newcommand{\M}[1]{\mathcal M_{{#1}}}
\newcommand{\Mm}[2]{\mathcal M_{{#1},{#2}}}
\newcommand{\Mb}[1]{\overline{\mathcal M}_{{#1}}}
\newcommand{\Mmb}[2]{\overline{\mathcal M}_{{#1},{#2}}}

\def\bb{\mathbb}

\def\phi{\varphi}

\def\Z{\bb{Z}}
\def\Q{\bb{Q}}
\def\R{\bb{R}}
\def\C{\bb{C}}
\def\A{\bb A_\C}

\def\co{\colon\thinspace}
\def\pu{\bullet}

\numberwithin{equation}{section}

\begin{document}

\author{Klaus Hulek}
\address{Leibniz Universit\"at Hannover, Institut f\"ur Algebraische Geometrie, Wel\-fen\-gar\-ten 1, D-30167 Hannover, Germany}
\email{hulek@math.uni-hannover.de}
\author{Orsola Tommasi}
\address{Leibniz Universit\"at Hannover, Institut f\"ur Algebraische Geometrie, Wel\-fen\-gar\-ten 1, D-30167 Hannover, Germany}
\email{tommasi@math.uni-hannover.de}

\title[Cohomology of the toroidal compactification of $\mathcal A_3$]{Cohomology of the toroidal compactification\\ of $\mathcal A_3$}
\date{March 29, 2009}

\begin{abstract}
We prove that the cohomology groups with rational coefficients of the Voronoi compactification $\AVOR 3$ of the moduli space of abelian threefolds coincide with the Chow groups of that space, as determined by Van der Geer.
\end{abstract}
\maketitle

\section{Introduction}

The moduli space $\ab g$ of principally polarized abelian varieties has several compactifications, notably the
Satake compactification $\ASAT g$ and various toroidal compactifications. Among the toroidal compactifications
the so called Voronoi compactification $\AVOR g$ is distinguished by the fact that it represents a geometrically
meaningful functor, as was shown by Alexeev \cite{A} and Olsson \cite{O}. Toroidal compactifications
are defined by suitable fans in the cone of semi-positive symmetric real $(g \times g)$ matrices and in the 
case of $\AVOR g$ the fan is given by the second Voronoi decomposition.
For a definition of the second Voronoi fan we refer the reader to \cite{V2b} or for a more
modern reference to \cite{AN}.
 A general discussion of toroidal compactifications 
of $\ab g$ can be found in the survey article \cite{HuSa}.

In genus $3$ all known toroidal compactifications of the moduli space $\ab 3$ of principally polarized abelian
varieties coincide with the Voronoi compactification $\AVOR 3$. 
We recall the explicit description of the 
second Voronoi decomposition in the case $g=3$ in Section \ref{rank3}.
A detailed description of the geometry of the space 
$\AVOR 3$ can be found in \cite{Tsushima}, and 
 the Chow ring of this space has been computed by
Van der Geer \cite{vdG}.

In this note we compute the cohomology groups with rational coefficients and prove that they
coincide with the Chow groups of this space. 
\begin{thm}\label{thm}
The Betti numbers of $\AVOR 3$ are  $b_0=b_{12}=1$, $b_2=b_{10}=2$, $b_4=b_8=4$ and $b_6=6$.
\end{thm}

Our approach is similar to that of \cite{vdG}, and is based on a study of the stratification of $\AVOR 3$ defined by the torus rank, which  we introduce in Section \ref{strat}.

We shall give the proof of the main result in Section \ref{proofmaintheo} 
modulo the computation of the cohomology of the various strata, which will be done in the subsequent sections.

As a corollary we obtain
\begin{cor}\label{cor}
The cycle map defines an isomorphism 
$$
CH^\pu(\AVOR 3)\otimes \Q \cong H^\pu(\AVOR 3;\Q)
$$
between the Chow ring and the cohomology ring of $\AVOR 3$ with rational coefficients.
\end{cor}
\begin{proof}
The Betti numbers coincide with the rank of the Chow groups as determined
by Van der Geer (\cite{vdG}). Since the intersection pairing is non-degenerate, the cycle map gives
an isomorphism.  
\end{proof}

Although this result is not particularly surprising, we could not find a reference to it in the literature, so we decided to fill the gap with this note.

We would like to remark that there are other possible approaches that yield the same result.
For instance, one can consider the Torelli map $\Mb 3\rightarrow \AVOR 3$ from the moduli space of Deligne--Mumford stable curves of genus $3$ to the toroidal compactification of $\ab 3$. 
The moduli space $\Mb 3$ has a stratification by topological type. Since the Torelli map for genus $3$ is surjective, we can stratify $\AVOR 3$ by taking the images of the strata of $\Mb 3$. It is easy to show that all strata of $\AVOR 3$ obtained in this way are isomorphic to finite quotients of products of moduli spaces $\Mm gn$ with $g\leq 3$ and $0\leq n\leq 2(g-3)$. Then one can use the known results about the cohomology of these spaces $\Mm gn$ to calculate the cohomology of $\AVOR 3$. 

In this note, we will work with the stack $\AVOR 3$ rather than the associated coarse moduli space. We recall
that $\AVOR 3$ is a smooth Deligne--Mumford stack. Hence the rational cohomology of the stack and the associated coarse
moduli space coincide.  

Finally we remark that the same techniques also apply to the (easier) case of genus $2$.

\begin{rem}
There is an isomorphism $$CH^\pu(\AVOR 2)\otimes\Q \cong H^\pu(\AVOR 2;\Q).$$
\end{rem}

\subsection*{Notation}\ \\
\begin{longtable}{lp{10cm}}
$\ab g$ & moduli stack of principally polarized abelian varieties of genus~$g$\\[1pt]
$\cX_g$ & universal family over $\ab g$\\[1pt]
$\ASAT g$ & Satake compactification of $\ab g$\\[1pt]
$\AVOR g$ &  Voronoi compactification of $\ab g$ \\[1pt]
$\XVOR g$& universal family over $\AVOR g$\\[1pt]
$\Mm gn$ & moduli stack of non-singular curves of genus $g$ with $n$ marked points\\[1pt]
\multicolumn{2}{l}{$\M g:=\Mm g0$}\\[1pt]
$\s_d$ & symmetric group in $d$ letters\\[1pt]
\end{longtable}
\smallskip

For every $g$, we denote by $\phi_g\co \AVOR g \rightarrow \ASAT g$ the natural map from the Voronoi to the Satake compactification.
Let $\pi_g: \XVOR g \to \AVOR g$ be the universal family, $q_g: \XVOR g \to \XVOR g/\pm 1$ the quotient
map from the universal family to the universal Kummer family and $k_g: \XVOR g/\pm 1 \to \AVOR g$
the universal Kummer morphism.

Throughout the paper, we work over the field $\C$ of complex numbers.

\subsection*{Acknowledgements} Partial support from DFG under grant Hu/6-1 is gratefully acknowledged.
We are also grateful to the referee for careful reading and valuable suggestions concerning the
presentation.

\section{Stratification and outline of the proof}\label{stratification}

\subsection{A stratification}\label{strat}
The object of this note is the rational cohomology of the toroidal compactification $\AVOR 3$ of the moduli space of abelian varieties of dimension $3$. 
We shall make use of a natural stratification of $\AVOR 3$ which was also used by Van der Geer \cite{vdG},
whose notation we adopt.

Recall that there is a natural map $\phi_3\co \AVOR 3 \rightarrow \ASAT 3$ to the Satake compactification. The moduli space $\ASAT 3$ admits a stratification $\ASAT 3 = \ab 3\sqcup \ab2\sqcup \ab1\sqcup\ab0$. This defines a filtration 
$\{\beta_t\}_{0\leq t\leq 3}$ on $\AVOR 3$, by setting 
$$\beta_t:=\phi_3^{-1}\left(\bigsqcup_{0\leq j\leq g-t}\ab j\right).$$
In other words, $\beta_t\subset\AVOR 3$ is the locus of semi-abelian varieties with torus rank at least $t$.

\subsection{Cohomology of the strata}\label{cohostrat}
We shall now state the results about the cohomology with compact support of the various strata. Proofs will be 
given in the subsequent 
Sections \ref{rank1} -- \ref{rank3}.

The stratum $\beta_0\setminus\beta_1$ of the filtration $\{\beta_t\}$ is $\ab 3$. Its cohomology was computed by Hain in \cite{Hain}. 

\begin{thm}\label{hain}
The rational cohomology groups with compact support of $\ab 3$ are given by
$$
H_c^k(\ab 3;\Q)=\left\{\begin{array}{ll}
\Q(-6) & k=12,\\
\Q(-5) & k=10,\\
\Q(-4) & k=8,\\
F & k=6,\\
0 & \text{otherwise},\\
\end{array}\right.
$$
where $F$ is a two-dimensional mixed Hodge structure which is an extension 
$$0\rightarrow \Q \rightarrow F\rightarrow \Q(-3)\rightarrow 0.$$
\end{thm}

\begin{proof}
This is a rephrasing of \cite[Thm~1]{Hain}, by using the isomorphism $$H_c^k(\ab 3;\Q)\cong H^{12-k}(\ab 3;\Q)^*\otimes \Q(-6)$$  given by Poincar\'e duality on the $6$-dimensional space $\ab 3$.
\end{proof}

The cohomology with compact support of the other strata is as follows.

\begin{prop}\label{propcohorank1}
The cohomology with compact support of $\beta_1\setminus\beta_2$ is given by 
$$
\begin{array}{l}
H_c^{10}({\beta_1\setminus\beta_2};\Q)=\Q(-5)\\
H_c^{8}({\beta_1\setminus\beta_2};\Q)=\Q(-4)^2\\
H_c^{6}({\beta_1\setminus\beta_2};\Q)=\Q(-3)^2\\
\end{array}
\ 
\begin{array}{l}
H_c^5({\beta_1\setminus\beta_2};\Q)=\Q \\
H_c^{4}({\beta_1\setminus\beta_2};\Q)=\Q(-2)\\
H_c^k({\beta_1\setminus\beta_2};\Q)=0 \text{    for }k\notin\{4,5,6,8,10\}.\\
\end{array}
$$
\end{prop}
For torus rank $2$ we obtain
\begin{prop}\label{propcohorank2}
The cohomology with compact support of the stratum $\beta_2\setminus\beta_3$ is given by
$$
\begin{array}{l}
H_c^{8}({\beta_2\setminus\beta_3};\Q)=\Q(-4)\\
H_c^{6}({\beta_2\setminus\beta_3};\Q)=\Q(-3)^2\\
H_c^{4}({\beta_2\setminus\beta_3};\Q)=\Q(-2)\\
\end{array}
\ \ \ \ 
\begin{array}{l}
H_c^2({\beta_2\setminus\beta_3};\Q)=\Q(-1) \\
H_c^k({\beta_2\setminus\beta_3};\Q)=0 \text{\ \ \ for }k\notin\{2,4,6,8\}.\\

\ \\
\end{array}
$$
\end{prop}

In the proof of the two propositions above, we make use of the following fact (see \cite{Tsushima}):
the natural map $\beta_1\rightarrow \ASAT 2$ factors through $\AVOR 2$, giving rise to the 
commutative diagram

\begin{equation*}
\xymatrix{
{\beta_1} \ar[r]^{k_2} \ar[rd]^{\pi_2}& {\AVOR 2}\\
&{\ASAT {2}}\ar@{<-}[u]_{\phi_2}
}
\end{equation*}
where $k_2\co \beta_1\cong(\XVOR2/\pm1)\rightarrow \AVOR2$ is the universal Kummer variety over $\AVOR 2$.

Finally, we use the toroidal description of $\AVOR 3$ to compute the cohomology of the stratum with torus rank $3$.  The corresponding result is
\begin{prop}\label{propcohorank3}
The cohomology groups $H_c^k(\beta_3;\Q)$ are trivial for degree $k$ different from $0,2,4,6$, and are given in the other cases by
$$
\begin{array}{llll}
H_c^6(\beta_3;\Q)=\Q(-3)&&&
H_c^2(\beta_3;\Q)=\Q(-1)
\\
H_c^4(\beta_3;\Q)=\Q(-2)^2&&&
H_c^0(\beta_3;\Q)=\Q.
\end{array}
$$
Moreover, the generators of these cohomology groups with compact support can be identified with the fundamental classes of the strata of $\AVOR 3$ corresponding to the cones $\sigma^{(3)}_{\text{local}}$, $\sigma^{(4)}_{I}$, $\sigma^{(4)}_{II}$, $\sigma^{(5)}$ and $\sigma^{(6)}$ (to be defined in Section~\ref{rank3}). 
\end{prop}

\subsection{Spectral sequences in cohomology}\label{spectral}

Our proofs of results on the cohomology of $\AVOR 3$ and its strata are based on an intensive use of long exact sequences and spectral sequences in cohomology with compact support. We shall recall the definition of the sequences we use most often in the proofs.

Since the cohomology with rational coefficients of a Deligne--Mumford stack coincides with that of its coarse moduli space,
in this section we will work with quasi-projective varieties. 
A more stack-theoretical approach can be obtained by recalling that $\AVOR g$ is the finite quotient of the fine moduli scheme $\AVOR g(n)$ of abelian varieties with level-$n$ structure for $n\geq3$. 
Then the same constructions can be obtained by working on $\AVOR g(n)$ equivariantly.

Recall that if $X$ is a quasi-projective variety and $Y$ a closed subvariety of $X$, then the inclusion $Y\hookrightarrow X$ induces a Gysin long exact sequence in cohomology with compact support:
\begin{equation*}
\cdots\rightarrow H_c^{k-1}(Y;\Q)\rightarrow H_c^{k}(X\setminus Y;\Q)\rightarrow H_c^k(X;\Q)\rightarrow H_c^{k}(Y;\Q)\rightarrow\cdots 
\end{equation*}

By functoriality of mixed Hodge structures  (\cite[Prop.~5.54]{PS},
this exact sequence respects mixed Hodge structures.
 
Next, assume we have  a filtration $\emptyset=Y_0\subset Y_1\subset Y_2\subset \cdots\subset Y_N = X$ by closed subvarieties of $X$. In this case, there is a spectral sequence  $E_r^{p,q}\Rightarrow H_c^{p+q}(X;\Q)$ associated to the filtration $\{Y_i\}$. The $E_1$ term is given by $E_1^{p,q}=H_c^{p+q}(Y_p\setminus Y_{p-1};\Q)$. This spectral sequence can be constructed by taking a compactification $\overline X$ of $X$ with border $S:=\overline X\setminus X$. Let us denote by  $\overline{Y_i}$ the closure of $Y_i$ in $\overline X$, and consider  the filtration $\{Y'_i:=\overline{Y_i}\cup S\}_{0\leq i\leq N}$ of the pair $(\overline X,S)$. In particular, one has $H^\pu(Y'_j,Y'_{j-1};\Q)=H^\pu_c(Y_j\setminus Y_{j-1};\Q)$ for all $j\geq 1$.  One can describe the spectral sequence associated to $\{Y_i\}$ as the spectral sequence associated to the bigraded exact couple $(D,E)$ with $D^{\alpha,\beta}=H^{\alpha+\beta}(\overline X,Y'_{\alpha-1};\Q)$ and $E^{\alpha,\beta}=H^{\alpha+\beta}(Y'_{\alpha},Y'_{\alpha-1};\Q)$, which converges to $H_c^{\alpha+\beta}(\overline X,S;\Q)=H^{\alpha+\beta}(X;\Q)$. Arguing as in \cite[Lemma~3.8]{arapura}, this ensures the compatibility with mixed Hodge structures by functoriality. For the definition of exact couples, see \cite[\S A.3.2]{PS}.

Note that the $d_1$ differentials of the spectral sequence in cohomology with compact support associated to $\{Y_i\}$ coincide with the differentials of the Gysin long exact sequences associated to the closed inclusions $Y_i\setminus Y_{i-1}\hookrightarrow Y_{i+1}\setminus Y_{i-1}$.

Leray spectral sequences play an intensive role in our computation of the cohomology of the strata $\beta_1$ and $\beta_2\setminus \beta_1$. Typically, we will be in  the following situation: let $X$ and $Y$ be quasi-projective varieties, and $f:\;X\rightarrow Y$ a fibration with fibres which are homotopy equivalent under proper maps to a fixed quasi-projective variety $B$. Let us denote by $\mathcal H^{(p)}$ the local system on $X$ induced by the $p$th cohomology group with compact support of the fibre of $f$.

In this situation, one can consider the Leray spectral sequence of cohomology with compact support associated to $f$. This is the spectral sequence $E_r^{p,q}\Rightarrow H_c^{p+q}(X;\Q)$ with $E_2^{p,q}\cong H^p(Y;\mathcal{H}^{(q)})$. Note that the Leray spectral sequence associated to $f$ respects Hodge mixed structures 
(e.g. see \cite[Cor.~6.7]{PS}).

\subsection{Proof of the main theorem}\label{proofmaintheo}
The results on the cohomology with compact support stated in Section~\ref{cohostrat} enable us to compute the cohomology of $\AVOR 3$ using the spectral sequence 
$E_\pu^{p,q}\Rightarrow H_c^{p+q}(\AVOR 3;\Q)$, $E_1^{p,q}=H_c^{p+q}(\beta_{3-p}\setminus\beta_{4-p};\Q)$
associated to the filtration $\beta_3\subset\beta_2\subset\beta_1\subset\beta_0=\AVOR 3$. 

\begin{lem}\label{pfthm1}
The $E_1$ term of the spectral sequence in cohomology with compact support associated to the filtration $\beta_3\subset\beta_2\subset\beta_1\subset\beta_0=\AVOR 3$ is  as given in Table~\ref{main_spseq}. The only non-trivial differential of this spectral sequence is $d_1^{2,3}\co E_1^{2,3}\rightarrow E_1^{3,3}$, which is injective. In particular, the spectral sequence degenerates at $E_2$.
\end{lem}

\begin{table}
\caption{\label{main_spseq} $E_1$ term of the spectral sequence converging 
to $H_c^\pu(\AVOR  3;\Q)=H^\pu(\AVOR  3;\Q)$}
$$
\begin{array}{r|ccccc}
q&&&&\\[6pt]
9&0&0&0&\Q(-6)\\
8&0&0&\Q(-5)&0\\
7&0&\Q(-4)&0&\Q(-5)\\
6&\Q(-3)&0&\Q(-4)^2&0\\
5&0&\Q(-3)^2&0&\Q(-4)\\
4&\Q(-2)^2&0&\Q(-3)^2&0\\
3&0&\Q(-2)&\Q&F\\
2&\Q(-1)&0&\Q(-2)&0\\
1&0&\Q(-1)&0&0\\
0&\Q&0&0&0\\
\hline
 &0&1&2&3&p\\
\left(\begin{smallmatrix}\text{torus}\\\text{rank}\end{smallmatrix}\right) 
& (3)&(2)&(1)&(0)
\end{array}
$$
\end{table}
\proof
The description of the $E_2$ term of the spectral sequence follows from the description of the compactly supported cohomology of the strata given in Section~\ref{cohostrat} and from the definition of the spectral sequence in Section~\ref{spectral}.

An inspection of the spectral sequence in Table~\ref{main_spseq} yields that $E_1^{p,q}$ (and hence $E_r^{p,q}$) is always trivial if $p+q$ is odd, with the exception of $E_1^{2,3}$ (hence possibly also $E_r^{2,3}$ for $r\geq 2$). Therefore, all differentials not involving $E_r^{2,3}$ terms are necessarily trivial, since they are maps either from or to $0$.  

This leaves us with only three possibly non-trivial differentials to investigate. The first two are the differentials $d_r^{2-r,2+r}\co E_r^{2-r,2+r}\rightarrow E_r^{2,3}$ for $r=1,2$. Note that in both cases, the Hodge structure on $E_r^{2,3}$ is pure of weight $0$, whereas the Hodge structure on $E_r^{2-r,2+r}$ is pure of weight $4$. Since the weights are different, the differential $d_r^{2-r,2+r}$ can only be the $0$ morphism.

Next, we investigate the differential $d_1^{2,3}\co E_1^{2,3}\rightarrow E_1^{3,3}$, which can have rank either $0$ or $1$.  Assume for the moment that $d_1^{2,3}$ is the $0$ morphism. Then the spectral sequence degenerates at $E_1$, so that $H^{5}(\AVOR 3;\Q)\cong H_c^{5}(\AVOR 3;\Q)=E_1^{2,3}=\Q$ holds. This means that the cohomology of $\AVOR 3$ in degree $5$ is pure of Hodge weight~$0$. 
But $\AVOR 3$ is a smooth proper stack, being the quotient by a finite group
of the stack $\AVOR 3(n)$ of principally polarized abelian varieties with a level-$n$ structure, which is represented by a smooth projective scheme for $n\geq 3$. In particular, the Hodge structure on $H^k(\AVOR3;\Q)$ is pure of weight $k$. 
Hence, the rank of $d_1^{2,3}$ must be $1$. Therefore, this differential is injective with cokernel isomorphic to $\Q(-3)$. This ensures $E_2^{2,3}=0$ and $E_2^{3,3}=\Q(-3)$. 
\qed 

Note that Lemma~\ref{pfthm1} directly implies that the cohomology of $\AVOR3$ is all algebraic, with Betti numbers as stated in Theorem~\ref{thm}.

In the remainder of this paper we will discuss the various strata defined by the torus rank and compute 
their cohomology.

\section{Torus rank $1$}\label{rank1}

To compute the cohomology with compact support of $\beta_1\setminus\beta_2$ we will use the map $k_2\co \beta_1\setminus\beta_2\rightarrow \ab 2$ realizing $\beta_1\setminus\beta_2$ as the universal Kummer variety over $\ab 2$. The fibre of $\beta_1\setminus\beta_2$ over a point parametrizing an abelian surface $S$ is $K:=S/\pm 1$. The cohomology of $K$ is one-dimensional in degree $0$ and $4$. The only other non-trivial cohomology group is $H^2(K;\Q)\cong \bigwedge^2H^1(S;\Q)$. 

To compute $H_c^\pu(\beta_1\setminus\beta_2;\Q)$, we consider the Leray spectral sequence associated to $k_2$. 
Note that the $0$th and the fourth cohomology group of the fibre induce trivial local systems on $\ab 2$.
Moreover, the second cohomology group of the fibre induces the rank $6$ local system $\bb V_{(1,1)}\oplus \Q(-1)$ on $\ab 2$. Here we denote by $\bb V_{(1,1)}$ the symplectic local system on $\ab 2$ determined by the irreducible representation of $\Sp(4,\Q)$ associated to the partition $(1,1)$. 

We start by determining the cohomology with compact support of $\ab2$ with values in the local system $\bb V_{(1,1)}$.

\begin{lem}\label{a2v11}
The rational cohomology groups with compact support of the moduli spaces $\M2$ and $\ab2$ with coefficients in $\bb V_{(1,1)}$ vanish in degree $k\neq 3$. In degree $3$, one has
$$H_c^3(\ab 2;\bb V_{(1,1)})=H_c^3(\M2;\bb V_{(1,1)})=\Q.$$
\end{lem}

\begin{proof}
We prove the claim about the cohomology of $\M2$ first.

Following the approach of \cite{G-2}, we use the forgetful map $p_2\co \Mm 22\rightarrow \M2$ to obtain information. Note that the fibre of $p_2$ is the configuration space of $2$ distinct points on a genus $2$ curve. The cohomology of $\Mm 22$, with the action of the symmetric group, was computed in \cite[Cor.~III.2.2]{OT-thesis}. This result allows us to conclude $H_c^3(\M2;\bb V_{(1,1)})=\Q$, $H_c^k(\M2;\bb V_{(1,1)})=0$ for $k\neq 3$. 
(Note that this is in agreement with the Hodge Euler characteristic of $\M2$ in the local system $\bb V_{(1,1)}$ computed in \cite[\S8.2]{G-2}.) 

Next, we determine $H_c^\pu(\ab2;\bb V_{(1,1)})$. To this end, we write $\ab 2$ as the disjoint union of the locus  $A_{1,1}$ of decomposable abelian surfaces, and the image of the Torelli map $t\co\M2\rightarrow \ab2$. 

Since the Torelli map is injective on the associated coarse moduli spaces, it induces an isomorphism between the cohomology of $\M2$ and that of $t(\M2)$ in every system of coefficients that is locally isomorphic to a $\Q$-vector bundle. 
Therefore, the Gysin long exact sequence with $\bb V_{(1,1)}$-coefficients associated to $A_{1,1}\hookrightarrow\ab 2$ yields
$$
H_c^{k-1}(A_{1,1};\bb V_{(1,1)})\rightarrow H_c^{k}(\M2;\bb V_{(1,1)})\rightarrow  H_c^k(\ab2;\bb V_{(1,1)})\rightarrow H_c^{k}(A_{1,1};\bb V_{(1,1)}).
$$
In Lemma~\ref{a11v11} below, we will show that $H_c^\pu(A_{1,1};\bb V_{(1,1)})$ is trivial. In view of the Gysin exact sequence above, this implies that $H_c^k(\ab2;\bb V_{(1,1)})$ is isomorphic to $H_c^{k}(\M2;\bb V_{(1,1)})$. This implies the claim.
\end{proof} 

\begin{rem}
Getzler's result would have been sufficient for the purposes of this note. This follows again from the fact that $\AVOR 3$ is a finite quotient of $\AVOR 3(n)$, so in particular its Hodge Euler characteristic determines the cohomology of the space as graded vector space with $\Q$-Hodge structures.
\end{rem}

\begin{lem}\label{a11v11}
The cohomology with compact support of $A_{1,1}$ in the local system of coefficients given by the restriction of $\bb V_{(1,1)}$ is trivial. 
\end{lem}
\begin{proof}
We consider the restriction of $k_2$ to $A_{1,1}$. Let $S=E_1\times E_2$ be an element of $A_{1,1}$, and let $K:=k_2^{-1}(S)$. Recall that $\bb V_{(1,1)}\oplus\Q(-1)$ is the local system $\mathcal H^{(2)}$ on $A_{1,1}$ induced by $H^2(K;\Q)$.
Therefore, the cohomology of $A_{1,1}$ with values in $\bb V_{(1,1)}\oplus \Q(-1)$ coincides with the cohomology of $A_{1,1}$ with values in the local system induced by the part of  $\bigwedge^2H^1(S;\Q)$ which is invariant under the symmetries of $E_1\times E_2$ and under the interchange of the two factors $E_1,E_2$ (which can be done topologically albeit not algebraically). Using the K\"unneth formula one sees that the latter local system is one-dimensional and induces the local system $\Q(-1)$. From this one obtains  $H_c^\pu(A_{1,1};\bb V_{(1,1)})=0$. 
\end{proof}

This allows us to show the following result, which directly implies that the cohomology with compact support of $\beta_1\setminus\beta_2$ is as stated in Proposition \ref{propcohorank1}.

\begin{proof}[Proof of Proposition~\ref{propcohorank1}.]
We compute the cohomology with compact support of $\beta_1\setminus\beta_2$ by using the Leray spectral sequence associated to the Kummer fibration $k_2\co\beta_1\setminus\beta_2\rightarrow\ab2$.

By the description of the fibre of $k_2$ given at the beginning of this section, the local systems $\mathcal H^{(0)}$ and $\mathcal H^{(4)}$ are the constant one, whereas $\mathcal H^{(2)}$ is the direct sum of the constant local system $\Q$ and $\bb V_{(1,1)}$.

The cohomology with compact support of $\ab 2$ is well known: it is one-dimen\-sion\-al in degree $4$ and $6$, and trivial elsewhere. This can be easily deduced from the results in \cite{Mu-curves} on the Chow ring of $\M2$. The cohomology of $\ab2$ in the local system $\bb V_{(1,1)}$ was computed in Lemma~\ref{a2v11}. From this, one obtains that
the $E_2$ term of the Leray spectral sequence in cohomology with compact support associated to $k_2$ is as in Table~\ref{specM2}.

From an inspection of the spectral sequence, one finds that all $E_2^{p,q}$ have pure Hodge structures, which have the same Hodge weight if and only if the sums $p+q$ coincide. Therefore, all differentials $d_r$ ($r\geq 2$) of the spectral sequence are morphisms between Hodge structures of different weight. Hence all differentials are trivial for this reason. This means that the spectral sequence degenerates at $E_2$, thus implying Proposition~\ref{propcohorank1}. 
\begin{table}
\caption{\label{specM2} $E_2$ term of the Leray spectral sequence converging to the cohomology with compact support of $\beta_1\setminus\beta_2$}
$$
\begin{array}{r|ccccc}
q&&&&\\[6pt]
4&0&\Q(-4)&0&\Q(-5)\\
3&0&0&0&0\\
2&\Q&\Q(-3)&0&\Q(-4)\\
1&0&0&0&0\\
0&0&\Q(-2)&0&\Q(-3)
\\\hline
&  3&4&5&6&p
\end{array}
$$
\end{table}
\end{proof}

\section{Torus rank $2$}\label{rank2}

Recall that $k_2\co\beta_1\rightarrow \AVOR 2$ is the universal family of Kummer varieties over $\AVOR 2$. Under this map, the elements of $\AVOR 3$ with torus rank $2$ are mapped to elements of $\AVOR 2$ of torus rank $1$. If we denote by $\beta'_t$ the stratum of $\AVOR 2$ of semi-abelian varieties of torus rank $\geq t$, we get a commutative diagram
$$\xymatrix@R=10pt@C=18pt{
{\AVOR 3}\ar@{<-^{)}}[d] &{\AVOR 2}\ar@{<-^{)}}[d] &{\AVOR 1}\ar@{<-^{)}}[d] \\
{\beta_2\setminus\beta_3} \ar[r]^{k_2} & {\beta'_1\setminus\beta'_2}\ar[r]^{k_1}& {\ab 1}\\
{}&{}&{}\\
{\pi_2^{-1}(\beta'_1\setminus\beta'_2)}\ar[uu]^{q_2}\ar[uur]_{\pi_2}&{{\mathcal X}_1}\ar[uu]^{q_1}\ar[uur]_{\pi_1}&{}
}$$

The map $\pi_2$ is the restriction of the universal family over $\AVOR 2$. In particular, the fibres of $\pi_2$ over points of $\beta'_1\setminus\beta'_2$ are rank $1$ degenerations of abelian surfaces, i.e. compactified $\C^*$-bundles over elliptic curves. 
A geometric description of these $\C^*$-bundles is given in \cite{Mu}.

We want to describe this situation in more detail. For this
consider the universal Poincar\'e bundle $\mathcal P\rightarrow \cX_1\times_{\ab 1}\hat{\cX_1}$ and let 
${\overline U}=\bb P(\mathcal P\oplus \mathcal{O}_{\cX_1\times_{\ab1}\cX_1})$ be the associated ${\bb P}^1$-bundle. 
Using the principal polarization we can naturally 
identify $\hat{\cX_1}$ and ${\cX_1}$, which we will do from now on.
We denote by 
$\Delta$ the union of the $0$-section and the $\infty$-section of this bundle. 
Set $U= {\overline U} \setminus \Delta$, which is simply the $\C^*$-bundle 
given by the universal Poincar\'e bundle $\mathcal P$ with the $0$-section removed and denote the 
bundle map by $f:U\rightarrow \cX_1\times_{\ab 1}\cX_1$.
Then there is a map $\rho: \overline U \to 
\beta_2\setminus\beta_3$ with finite fibres. 
Note that the two components of $\Delta$ are identified under the map $\rho$. The restriction of $\rho$ to both $U$ and to
$\Delta$ is given by a finite group action, although the group is not the same in the two cases (see the discussion below). 

We now consider the situation over a fixed point $[E] \in \ab1$.
For a fixed degree $0$ line bundle ${\mathcal L}_0$ on $E$ 
the preimage $f^{-1}(E \times \{{\mathcal L}_0 \})$ is a semi-abelian surface, namely 
the $\C^*$-bundle given by 
the extension corresponding to  ${\mathcal L}_0 \in \hat{E}$.  This semi-abelian surface admits a Kummer involution
$\iota$ which acts as $x \mapsto -x$ on the base $E$ and by $t \mapsto 1/t$ on the fibre over the origin. The Kummer involution $\iota$ is defined universally on $U$.

Consider the two involutions $i_1, i_2$ on $\cX_1\times_{\ab1}\cX_1$ defined by $i_1(E,p,q)=(E,-p,-q)$ and $i_2(E,p,q)=(E,q,p)$ for every elliptic curve $E$ and every $p,q\in E$. These two involutions lift to involutions $j_1$ and $j_2$ on $U$ that act trivially on the fibre of $f\co U\rightarrow \cX_1\times_{\ab1}\cX_1$ over the origin.

\begin{lem}
The diagram
\begin{equation}\label{poincare}
\xymatrix{
U \ar[r]\ar[d]^{\rho|_U} & {\cX_1\times_{\ab1}\cX_1}\ar[d]^{\rho'}\\
{(\beta_2\setminus\beta_3) \setminus \rho(\Delta)} \ar[r] & {\Sym^2_{\ab1}(\cX_1/\pm 1)},\\
}
\end{equation}
where $\rho':{\cX_1\times_{\ab1}\cX_1}\to{\Sym^2_{\ab1}(\cX_1/\pm 1)}$ is the natural map, is commutative. Moreover 
$\rho|_U\co U\rightarrow \rho(U)\subset\beta_2\setminus\beta_3$ is the quotient of $U$ by the subgroup of the automorphism group of $U$ generated by $\iota, j_1$ and $j_2$.
\end{lem}

\begin{proof}
Since the map $\rho'$ in the diagram~\eqref{poincare} has degree $8$ and $\iota, j_1,j_2$ generate a subgroup of order $8$ of the automorphism group of $U$, it suffices to show that the map $\rho|_U$ factors through each of the involutions $\iota$ and $j_1, j_2$. 
 
Recall that the elements of $\beta_2\setminus\beta_3$ correspond to rank $2$ degenerations of abelian threefolds. 
More precisely, every point of $\rho(U)$ corresponds to a degenerate abelian threefold $X$ whose 
normalization is a ${\bb P}^1 \times {\bb P}^1$-bundle, namely the 
compactification of a product of two $\C^*$-bundles on the elliptic curve $E$ given by $k_1\circ k_2([X])$. 
The degenerate abelian threefold itself is given by identifying 
the $0$-sections and the $\infty$-sections of the ${\bb P}^1 \times {\bb P}^1$-bundle. This identification 
is determined by a complex parameter, namely the point on a fibre of $U\rightarrow \cX_1\times_{\ab 1}\cX_1$.

Since a degree $0$ line bundle ${\mathcal L}_0$ and its inverse define isomorphic semi-abelian surfaces and since
the role of the two line bundles is symmetric, the map $\rho|_U$ factors through $\iota$ and$j_2$. 
Since $j_1$ is the commutator of $\iota$ and $j_2$ the map $\rho|_U$ also factors through $j_1$.
\end{proof}

A consequence of the lemma above is that the cohomology with compact support of $\rho(U)$ can be computed by taking the invariant part of the cohomology of the total space of the $\C^*$-bundle $f\co U\rightarrow \cX_1\times_{\ab1}\cX_1$. 
Hence, the invariant part of the Leray spectral sequence associated to $f$
gives a Leray spectral sequence converging to $H_c^\pu(\rho(U);\Q)$.
Thus, we have to consider the part of $E_2^{p,q}(f)=H_c^q(\C^*;\Q)\otimes H_c^p(\cX_1\times_{\ab1}\cX_1;\Q)$ that is invariant under the action of $\iota, j_1$ and $j_2$.

Since $j_1$ and $j_2$ both fix the fibre of $f$ over the origin, they act trivially on the cohomology of $\C^*$. Instead, the Kummer involution $\iota$ acts as the identity on $H_c^2(\C^*;\Q)$ and as the alternating representation on $H_c^1(\C^*;\Q)$.

The action of $\iota$, $j_1$ and $j_2$ can be determined by considering the induced actions on $\cX_1\times_{\ab1}\cX_1$. Here one uses that all three involutions respect the map $\cX_1\times_{\ab1}\cX_1\rightarrow\ab1$, whose fibre over $[E]\in\ab1$ is isomorphic to $E\times E$. 
Note in particular that the involution $(E,p,q)\leftrightarrow (E,-p,q)$ induced by $\iota$ acts as the alternating representation on the linear subspace $\bigwedge^2H_c^1(E;\Q)\subset H_c^2(E\times E;\Q)$, on which $i_1$ and $i_2$ both act trivially.

This discussion yields that the invariant part of the spectral sequence $E_2$ term is as shown in Table~\ref{t:rank2}.
 
\begin{table}\caption{\label{t:rank2} $E_2$ term of the spectral sequence converging to the cohomology with compact support of $\rho(U)$} 
$$
\begin{array}{r|cccccc}
q&&&&\\[6pt]
2&\Q(-2)&0&\Q(-3)&0&\Q(-4)\\
1&0&0&\Q(-2)&0&0
\\\hline
& 2&3&4&5&6&p
\end{array}
$$
\end{table}

\begin{lem}
The cohomology groups with compact support of $\rho(U)$ are $1$-dimen\-sion\-al in degree $6$ and $8$ and trivial otherwise.
\end{lem}

\begin{proof}
It suffices to show that the differential $d_2^{2,2}\co E_2^{2,2}\rightarrow E_2^{4,1}$ in Table~\ref{t:rank2} is an isomorphism. 

To describe the differential $d_2^{2,2}$ geometrically, it is useful to consider the restriction of the Torelli map $t\co\Mb 3\rightarrow \AVOR 3$ to the preimage of $\rho(U)$. Moreover, one can use the stratification of $\Mb3$ by topological type to describe $\beta_2$ and $\rho(U)$. In particular, this allows one to find a geometric generator for $H_c^4(\rho(U);\Q)$. 

Consider stable curves $C_1\cup C_2\cup C_3$, where the component $C_1$ is smooth of genus $1$, the component $C_2$ is a smooth rational curve and the component $C_3$ is a rational curve with exactly one node, satisfying $\#(C_1\cap C_2)=1$, $\#(C_1\cap C_3)=0$ and $\#(C_2\cap C_3)=2$.

Denote by $G$ the closure in $t^{-1}(\rho(U))$ of the locus of such curves, and denote by $t_*[G]$ the push-forward to $\rho(U)$ of the cycle class of $G$. Then the fundamental class of $t_*[G]$ generates $H_c^4(\rho(U);\Q)$.

Recall that the locus in $\Mb 3$ of irreducible curves with two nodes maps surjectively to $\beta_2$ under the Torelli map.
Moreover, all curves in $\Mb 3$ that have two nodes and map to $\beta_2$
 can be constructed by taking a stable curve of genus $1$ with $4$ marked points and identifying the marked points pairwise.
There is a well known relation between cycle classes of dimension $2$ in $\Mmb 14$, called Getzler's relation (see \cite{G-rela}). This relation is $\s_4$-invariant and it induces a relation between dimension $2$ cycles in $t^{-1}(\beta_2)$, which if pushed forward under $t$ induces a relation in $H_c^4(\beta_2;\Q)$. The latter relation involves non-trivially the push-forward of the fundamental class of $\overline G\subset t^{-1}(\beta_2)$. In particular, restricting to $\rho(U)\subset\beta_2$ yields that $t_*[G]$ vanishes in $H_c^4(\rho(U);\Q)$. Hence, the differential $d_2^{2,2}$ must be an isomorphism.
\end{proof}

\begin{rem}
There is also another way to see that the differential $d_2^{2,2}\co E_2^{2,2}\rightarrow E_2^{4,1}$ in Table~\ref{t:rank2} is an isomorphism.
Namely, one can compactify the $\C^*$-bundle $U$ to the ${\bb P^1}$-bundle ${\overline U}=\bb P(\mathcal P\oplus \mathcal{O}_{\cX_1\times_{\ab1}\cX_1})$ and compute the invariant part of the exact sequence in rational cohomology
of the pair $({\overline U},\Delta)$. This then shows that the invariant part of $H^4_c(U;\Q)$ vanishes as claimed.
We decided to include the above proof involving Getzler's relation since the relation to $\Mb3$ is of independent interest. 
\end{rem}

\begin{proof}[Proof of Proposition~\ref{propcohorank2}]
We compute the cohomology with compact support of $\beta_2\setminus\beta_3$ by exploiting the Gysin long exact sequence associated to the inclusion $\rho(\Delta)\hookrightarrow(\beta_2\setminus\beta_3)$:
\begin{equation}\label{gysin2}
\cdots\rightarrow H_c^{k-1}(\rho(\Delta);\Q)\rightarrow H_c^{k}(\rho(U);\Q)\rightarrow H_c^k(\beta_2\setminus\beta_3;\Q)\rightarrow H_c^{k}(\rho(\Delta);\Q)\rightarrow\cdots 
\end{equation}

The map $\rho$ identifies the two components of $\Delta$, each of which is isomorphic to $\cX_1 \times_{\ab1} \cX_1$. Moreover,
it factors through the finite group $G$ generated by the following three involutions: the involution which interchanges
the two factors of $\cX_1 \times_{\ab1} \cX_1$, the involution which acts by $(x,y)\mapsto (-x,-y)$ on each fibre $E \times E$ and finally
the involution which acts by $(x,y) \mapsto (x+y,-y)$. This can be read off from the
construction of the toroidal compactification (see \cite[Section I]{HuSa} for an outline of this
construction. Also note that the stratum $\Delta$ corresponds to
the stratum in the partial compactification in the direction of the $1$-dimensional cusp associated to a
maximal-dimensional cone in the second Voronoi decomposition for $g=2$. A detailed
description can be found in \cite[Part I, Chapter 3]{HKW}).

Hence
$$H^\pu_c(\rho(\Delta);\Q)\cong H^\pu_c(E \times E/G;\Q) \otimes H^\pu_c(\A^1;\Q).$$
A straightforward calculation shows that the $G$-invariant cohomology of $E\times E$ has rank $1$ in even dimension and vanishes otherwise.
In particular this quotient behaves cohomologically like ${\mathbb P}^2$.

Since $H_c^k(\rho(U);\Q)$ and $H_c^k(\rho(\Delta);\Q)$ both vanish if $k$ is odd, the exact sequence \eqref{gysin2} splits into short exact sequences
$$0\rightarrow H_c^{k}(\rho(U);\Q)\rightarrow H_c^k(\beta_2\setminus\beta_3;\Q)\rightarrow H_c^{k}(\rho(\Delta);\Q)\rightarrow 0.$$
This implies the claim.
\end{proof}

\begin{rem}
We would like to take this opportunity to correct a slight error in \cite[3.8]{vdG} where it was claimed that the map $\rho$ factors through
$\Sym^2_{\ab1}(\cX_1/\pm 1)$ rather than through the quotient by $G$.
This, however, does not effect the results of \cite{vdG}.
\end{rem}

\section{Torus rank $3$}\label{rank3}

The stratum $\beta_3\subset \AVOR 3$ lying over $\ab 0\subset\ASAT 3$ is entirely determined by the fan of the toroidal compactification. For this we first have to describe the Voronoi fan $\Sigma$ in genus $3$. 

Consider the free abelian group $\bb L_3\cong \Z^3$ with generators $x_1,x_2,x_3$ and let $\bb M_3=\Sym_2(\bb L_3)$. 
Then $\bb M_3$ is isomorphic to the space of $3\times 3$ integer symmetric matrices with respect to the 
basis $x_i$ via the map which assigns to a matrix $A$ the quadratic form $^{t}xAx$.
We shall use the basis of $\bb M_3$ given by the forms $U_{i,j}^*$, $1\leq i\leq j\leq 3$ given by 
$$U_{i,j}^*=2^{\delta_{i,j}}x_ix_j.$$
Let $\Sym_2^{\geq 0}(\bb L_3\otimes\R)$ be the cone of positive semidefinite forms in $\bb M_3\otimes\Q$.
The group $\GL(3,\Z)$ acts on $\Sym_2^{\geq 0}(\bb L_3\otimes\R)$ by
$$\GL(3,\Z)\ni g:\ M\longmapsto {}^tg^{-1}Mg^{-1}.$$
Let
$$
\sigma^{(6)}:=\R_{\geq 0}\alpha_1+\R_{\geq 0}\alpha_2+\R_{\geq 0}\alpha_3+\R_{\geq 0}\polybeta_1+\R_{\geq 0}\polybeta_2+\R_{\geq 0}\polybeta_3,
$$
where $\alpha_i=x_i^2$ for all $i=1,2,3$ and $\polybeta_i=(x_j-x_k)^2$ for $\{i,j,k\}=\{1,2,3\}$.
Since the forms $\alpha_j,\polybeta_i$ form a basis of $\bb M_3$, this is a basic $6$-dimensional cone in 
$\Sym_2^{\geq 0}(\bb L_3\otimes\R)$. 

The Voronoi fan in genus $3$  is the fan $\Sigma$ in $\Sym_2^{\geq 0}(\bb L_3\otimes\R)$
given by $\sigma^{(6)}$ and all its faces, together with their $\GL(3,\Z)$-translates. We use the notation
$$\sigma^{(6)}=\alpha_1*\alpha_2*\alpha_3*\polybeta_1*\polybeta_2*\polybeta_3,$$
and similarly for the faces of $\sigma^{(6)}$. 

To describe $\AVOR 3$, we have to know all possible $\GL(3,\Z)$-orbits of  $\sigma^{(6)}$ and its faces. An $i$-dimensional cone corresponds to a $(6-i)$-dimensional stratum in $\AVOR 3$. Since strata of dimension at least $4$ necessarily lie over $\ab l$ with $l\geq 1$, we only need to know the orbits of cones of dimension $\leq 3$.

The following lemma can be proved using the methods of \cite{Tsushima} (see \cite[Chapter~3]{Erdenberger}).

\begin{lem}\label{class}
There are two $\GL(3,\Z)$-orbits of $3$-dimensional cones, represented by the cones 
$$
\sigma^{(3)}_{\text{local}}=\alpha_1*\alpha_2*\alpha_3,\ \  \
\sigma^{(3)}_{\text{global}}=\alpha_1*\alpha_2*\polybeta_3.
$$
The stratum associated to $\sigma^{(3)}_{\text{local}}$ lies over $\ab0$, that associated to $\sigma^{(3)}_{\text{global}}$ lies over $\ab 1$.

There are two $\GL(3,\Z)$-orbits of $4$-dimensional cones, given by
$$
\sigma^{(4)}_{I}=\alpha_1*\alpha_2*\alpha_3*\polybeta_1,\ \ \ 
\sigma^{(4)}_{II}=\alpha_1*\alpha_2*\polybeta_1*\polybeta_2.
$$
In dimension $5$ and $6$ there is only one $\GL(3,\Z)$-orbit. The strata of all cones of dimension at least $4$
lie over $\ab0$.
\end{lem}

Let
$$\bb H_3 = \{\tau=(\tau_{i,j})_{1\leq i,j\leq 3}: \tau={}^t\tau, \im\tau>0\}$$
be the Siegel upper half plane of genus $3$. We consider the rank $6$ torus $T=T^6$ with coordinates 
$$t_{i,j}=e^{2\pi \sqrt{-1} \tau_{i,j}}\ (1\leq i,j\leq 3).$$

These coordinates correspond to the dual basis of the basis $U_{i,j}^*$. 
If $\sigma^{(l)}$ is an $l$-dimensional cone in $\Sigma$ then, since the fan $\Sigma$ is basic, it follows that
the associated affine variety $T_{\sigma^{(l)}} \cong \C^l\times(\C^*)^{6-l}$. The corresponding stratum in $\AVOR 3$ is then a quotient of $\{(0,0,0)\}\times(\C^*)^{6-l}$ by a finite group.
We consider the torus embedding $T\hookrightarrow T_{\sigma^{(6)}}\cong \C^6$, where the latter isomorphism holds since $\sigma^{(6)}$ is a basic cone of dimension $6$. Let $T_1,\dots,T_6$ be the coordinates of $\C^6$ corresponding to the basis $\alpha_1,\dots,\polybeta_3$. 
If one computes the dual basis of $\alpha_1,\dots,\polybeta_3$ in terms of the dual basis of $U_{i,j}^*$, one obtains that the torus embedding $T\hookrightarrow \C^6$ is given by
$$ 
\begin{array}{lllll}
T_1=t_{1,1}t_{1,3}t_{1,2},&&
T_2=t_{2,2}t_{2,3}t_{1,2},&&
T_3=t_{3,3}t_{1,3}t_{2,3},\\
T_4=t_{2,3}^{-1},&&
T_5=t_{1,3}^{-1},&&
T_6=t_{1,2}^{-1}.
\end{array}
$$

Let us start by considering the stratum associated to 
$$\sigma^{(3)}_{\text{local}}=\alpha_1*\alpha_2*\alpha_3.$$

Let $S_1,S_2$ and $S_3$ be coordinates corresponding to $\alpha_1,\alpha_2$ and $\alpha_3$, and let $t_{2,3}, t_{1,3},$ $t_{1,2}$ be as above. Then 
$$ T_{\sigma^{(3)}_{\text{local}}}\cong \C^3\times(\C^*)^3\subset \C^6=T_{\sigma^{(6)}}$$
with coordinates $S_1,S_2,S_3,t_{2,3}^{-1}, t_{1,3}^{-1}, t_{1,2}^{-1}$,
where the inclusion is defined by considering $\sigma^{(3)}_{\text{local}}$ as a face of $\sigma^{(6)}$. 

The stratum which we add is $\{(0,0,0)\}\times(\C^*)^3$ modulo a finite group $G=G_{\sigma^{(3)}_{\text{local}}}$, namely the 
stabilizer of the cone $\sigma^{(3)}_{\text{local}}$ in $\GL(3,\Z)$.
In order to understand the action of the group $G$ explicitly, we recall that it is naturally a subgroup of
the parabolic subgroup which belongs to the standard $0$-dimensional cusp 
$$P=\left\{\left(\begin{array}{c|c}\;g&0\\[6pt]\hline\\[-6pt] 0&{}^tg^{-1}\end{array}\right): g\in\GL(3,\Z)\right\} \cong \GL(3,\Z) \subset \Sp(6,\Z).$$ 

\begin{lem}\label{s3loc}
The stratum associated to $\sigma^{(3)}_{\text{local}}$ is an affine variety $Y_{\text{local}}^{(3)}=(\C^*)^3/G$ whose only non-trivial cohomology with compact support is in degree $6$.
\end{lem}

\begin{proof}
Since the stratum associated to $\sigma^{(3)}_{\text{local}}$ is the quotient of the smooth variety $(\C^*)^3$ by a finite group, its cohomology and cohomology with compact support are related by Poincar\'e duality. Hence, it suffices to show that the rational cohomology of the stratum is concentrated in degree $0$.

Denote by $T^3$ the rank $3$ torus with coordinates  $(v_1,v_2,v_3)=(t_{2,3}^{-1}, t_{1,3}^{-1}, t_{1,2}^{-1})$.
The stratum which we add for $\sigma^{(3)}_{\text{local}}$ is then isomorphic to $T^3/G$. Since $\sigma^{(3)}_{\text{local}}=\alpha_1*\alpha_2*\alpha_3$ with $\alpha_i=x_i^2$, we see that the group $G$ 
is the group generated by the permutations of the $x_i$ and the involutions $(x_1,x_2,x_3)\mapsto(\epsilon_1x_1,\epsilon_2x_2,\epsilon_3x_3)$ with $\epsilon_i=\pm1$.
Note that the element $-id$ acts trivially both on $\bb H_3$ and on $\bb M_3$. Hence the group $G$ is an extension
$$1\rightarrow (\Z/2\Z)^2\rightarrow G \rightarrow \s_3 \rightarrow 1,$$
where $\s_3$ denotes the symmetric group in $3$ letters. Next, we have to analyze how this group acts on $\bb H_3$ and on the torus $T^3$. The permutation of $x_i$ and $x_j$ interchanges $\tau_{i,k}$ and $\tau_{j,k}$ but fixes $\tau_{i,j}$. Hence $\s_3$ also acts as group of permutations on the coordinates of $T^3$. The action of the involutions generating $(\Z/2\Z)^2$ can be seen for example from
$$
\begin{pmatrix}
-1&0&0\\
0&1&0\\
0&0&1
\end{pmatrix}
\begin{pmatrix}
\tau_{1,1}&\tau_{1,2}&\tau_{1,3}\\
\tau_{1,2}&\tau_{2,2}&\tau_{2,3}\\
\tau_{1,3}&\tau_{2,3}&\tau_{3,3}
\end{pmatrix}
\begin{pmatrix}
-1&0&0\\
0&1&0\\
0&0&1
\end{pmatrix}
=
\begin{pmatrix}
\tau_{1,1}&-\tau_{1,2}&-\tau_{1,3}\\
-\tau_{1,2}&\tau_{2,2}&\tau_{2,3}\\
-\tau_{1,3}&\tau_{2,3}&\tau_{3,3}
\end{pmatrix}.
$$

Hence, the involution $(x_1,x_2,x_3)\leftrightarrow (x_1,x_2,-x_3)$ induces the involution given by $(v_1,v_2,v_3)\leftrightarrow (v_1^{-1},v_2^{-1},v_3)$ and similarly for the other involutions. 
This allows us to describe the quotient $T^3/G$ explicitly, as given by the image of the map
$$
\begin{array}{c@{\;\;}c@{\;\;}l}T ^3 \cong (\C^*)^3&\longrightarrow&\C^4\\
(v_1,v_2,v_3)&\longmapsto & (u_1+u_2+u_3,u_1u_2+u_1u_3+u_2u_3,u_1u_2u_3,u_4)=(s_1,s_2,s_3,t),
\end{array}$$
where 
$$u_1=v_1+\frac1{v_1},\ \ 
u_2=v_2+\frac1{v_2},\ \ 
u_3=v_3+\frac1{v_3},\ \ 
u_4=\left(v_1-\frac1{v_1}\right)\left(v_2-\frac1{v_2}\right)\left(v_3-\frac1{v_3}\right).
$$

Then the image is the hypersurface $W\subset\C^4$ given by
$$\frac{t^2}{4}-(\frac{s_3}{2}+2s_1)^2+(s_2+4)^2=0. $$

Note that $W$ is a cone with vertex the line $t=\frac{s_3}{2}+2s_1=s_2+4=0$ in $\C^4$ over a plane projective conic. Then the claim follows from the contractibility of $W$.

Alternatively, one can also show that the cohomology $H^\pu(T^3/G;\Q)$ is concentrated in degree $0$, by proving that the only cohomology in $H^\pu(T^3;\Q)$ which is fixed under the group $G$ is in degree $0$. 
\end{proof}

The situation with the lower-dimensional strata is similar:

\begin{lem}\label{onlyh0}
Let $\sigma^{(l)}$ be an $l$-dimensional subcone of $\alpha_1*\alpha_2*\alpha_3*\polybeta_1*\polybeta_2*\polybeta_3$, with $l\geq 4$. Then the stratum of $\polybeta_3$ associated to $\sigma^{(l)}$ has non-trivial cohomology with compact support only in the maximal degree $2(6-l)$. 
\end{lem}

\begin{proof}
Recall that all $\GL(3,\Z)$-orbits of $\sigma^{(l)}$ were described in Lemma~\ref{class}. Hence it suffices to consider the cases in which $\sigma^{(l)}$ is one of the following cones:  
$\sigma^{(4)}_{I}$, $\sigma^{(4)}_{II}$, $\sigma^{(5)}:=\alpha_1*\alpha_2*\alpha_3*\polybeta_1*\polybeta_2$ and $\sigma^{(6)}$.

As mentioned above, if $\sigma^{(l)}$ is an $l$-dimensional cone in $\Sigma$ then we have $T_{\sigma^{(l)}}=\C^l\times(\C^*)^{6-l}$, because the fan $\Sigma$ is basic. The corresponding stratum in $\AVOR 3$ is then a quotient of $\{(0,0,0)\}\times(\C^*)^{6-l}\cong (\C^*)^{6-l}$ by a finite group $G_{\sigma^{(l)}}$. To prove the claim, it suffices to show that the part of the cohomology of $(\C^*)^{6-l}$ which is invariant for the action of $G$ coincides with $H^0((\C^*)^{6-l};\Q)$. 
Since $(\C^*)^{6-l}$ is smooth, the result about cohomology with compact support will follow from Poincar\'e duality.

For instance, consider the case of $\sigma^{(4)}_{II}$. Using toric coordinates, one finds that the corresponding stratum is given by a quotient of $(\C^*)^2$ by the action of the finite group $\Z/2\Z\times\s_3$. The factor $\s_3$ acts on $\polybeta_1*\alpha_2*\alpha_3$ by permuting $\polybeta_1$, $\alpha_2$ and $\alpha_3$, whereas the action of the factor $\Z/2\Z$ is generated by the involution $x_1\leftrightarrow -x_1$. One can compute explicitly the action of $\Z/2\Z\times\s_3$ and prove $(H^\pu((\C^*)^2;\Q))^{\Z/2\Z\times\s_3}=H^0((\C^*)^2;\Q)$.

Analogous considerations yield the claim in the case of the other strata. 
\end{proof}

Concluding, the proof of Proposition \ref{propcohorank3} now follows 
from Lemmas~\ref{s3loc} and \ref{onlyh0}.

\bibliographystyle{amsalpha}

\begin{thebibliography}{vdG}

\bibitem[A]{A}V.~Alexeev, {\it Complete moduli in the presence of
    semiabelian group action.}  Ann. of Math. (2)  {\bf 155}  (2002), 611--708.
\bibitem[AN]{AN}V.~Alexeev, I.~Nakamura, {\it On Mumford's construction of degenerating abelian varieties.}
 Tohoku Math. J. (2) {\bf 51}  (1999), 399--420.  
\bibitem[Ar]{arapura} D.~Arapura, {\it The Leray spectral sequence is motivic}.  Invent. Math.  {\bf 160}  (2005),  no. 3, 567--589.
\bibitem[E]{Erdenberger} C.~Erdenberger, {\it A finiteness result for Siegel modular threefolds}. Ph. D. Thesis, Leibniz Universit\"at Hannover (2007). Available at \texttt{http://www.iag.uni-hannover.de/$\sim$ag-iag/} \texttt{data/phdthesis\_erdenberger.pdf}

\bibitem[vdG]{vdG} G.~van der Geer, {\it The Chow ring of the moduli 
space of abelian threefolds}, J. Algebraic Geom. {\bf 7} (1998), 753--770.
\bibitem[G1]{G-rela}
E.~Getzler,
{\it Intersection theory on $\overline{\mathcal M}_{1,4}$ and elliptic Gromov--Witten invariants}
 J. Amer. Math. Soc.  {\bf 10}  (1997),  no. 4, 973--998.
\bibitem[G2]{G-2}
E.~Getzler,
{\it Topological recursion relations in genus {$2$}},
in: {Integrable systems and algebraic geometry (Kobe/Kyoto, 1997)},
{World Sci. Publ., River Edge, NJ},
{1998},
{73--106}.
\bibitem[H]{Hain}
R.~Hain,
{\it The rational cohomology ring of the moduli space of abelian $3$-folds},
Math. Res. Lett.  {\bf 9}  (2002),  no. 4, 473--491.
\bibitem[HKW]{HKW} K.~Hulek, C.~Kahn, S.~H.~Weintraub,
{\it Moduli spaces of abelian surfaces: compactification, degenerations, and theta functions}, de Gruyter Expositions in Mathematics
{\bf 12}. Walter de Gruyter \& Co., Berlin, 1993. 
\bibitem[HS]{HuSa} K.~Hulek, G.~K.~Sankaran, 
{\it The geometry of Siegel modular varieties.} 
Higher dimensional birational geometry (Kyoto, 1997),  Adv. Stud. 
Pure Math., {\bf 35}, 89--156, Math. Soc. Japan, Tokyo, 2002.
\bibitem[M1]{Mu} D.~Mumford, {\it On the Kodaira dimension of the Siegel
modular variety}. In: Algebraic geometry -- open problems, proceedings,
Ravello 1982, eds. C.~Ciliberto, F.~Ghione and F.~Orecchia. Lecture
Notes in Mathematics {\bf 997}, Springer-Verlag, Berlin-New
York, 1983, 348--375.
\bibitem[M2]{Mu-curves} D.~Mumford, {\it Towards an enumerative geometry of the moduli space of curves}. 
In:  Arithmetic and geometry, Vol. II.
Progr. Math., {\bf 36}, Birkh\"auser Boston, Boston, MA, 1983,  271--328.
\bibitem[O]{O} M.~Olsson, {\it Compactifying moduli spaces for abelian varieties.} 
Lecture Notes in Mathematics, {\bf 1958}. Springer-Verlag, Berlin, 2008. 
\bibitem[PS]{PS} C.~A.~M.~Peters, J.~H.~M.~Steenbrink, \emph{Mixed Hodge Structures}.
Ergebnisse der Mathematik und ihrer Grenzgebiete. 3. Folge,
\textbf{52},
Springer-Verlag,
Berlin, 2008.


\bibitem[T]{OT-thesis}
O.~Tommasi,
{\it Geometry of discriminants and cohomology of moduli spaces}.
Ph.D. thesis, Radboud University Nijmegen (2005).
\\
Available at 
  \texttt{http://webdoc.ubn.ru.nl/mono/t/tommasi\_o/geomofdia.pdf}
\bibitem[Ts]{Tsushima} R.~Tsushima, {\it  A formula for the dimension of spaces of Siegel cusp forms of degree three}, Amer. J. Math.  {\bf 102}  (1980), no. 5, 937--977.
\bibitem[V]{V2b} G.~F.~Voronoi, {\it Nouvelles applications des param\`etres
continus \`a la th\'eorie des formes quadratiques. Deuxi\`eme
m\'emoire. Recherches sur les parall\'elo\`edres primitifs. Seconde
partie. Domaines de formes quadratiques correspondant aux diff\'erents
types de parall\'elo\`edres primitifs},
J. Reine Angew. Math. {\bf 136} (1909), 67--178.

\end{thebibliography}

\end{document}